\documentclass[12pt,reqno]{amsart}
\usepackage{amsmath,amsthm,amsbsy,amssymb}
\advance \topmargin by-\headheight
\advance \topmargin by-\headsep
\evensidemargin \oddsidemargin
\theoremstyle{plain}

\def\bb{\mathbf{b}}
\def\gg{\mathbf{g}}

\newtheorem{thm}{Theorem}
\newtheorem{lem}[thm]{Lemma}

\theoremstyle{definition}

\theoremstyle{remark}

\theoremstyle{plain}

\newcommand{\norm}[1]{\left|\hspace*{2.5pt} \!\!\left| #1\right|\hspace*{2.5pt}
\!\!\right|}

\numberwithin{equation}{section}

\begin{document}
\baselineskip 18pt

\title[On a Problem of Mordell with Primitive Roots]{On a Problem of Mordell
with Primitive Roots}
\author{Cristian Cobeli}

\address{Cristian Cobeli, 
Mathematics Research Institute of the Romanian Academy,
P.O. Box 1-764, Bucharest, 70700, Romania}

\email{cristian.cobeli@imar.ro}
\date{\today}

\maketitle
\begin{abstract}
   We  consider the sums of the form
\begin{equation*}
S=\sum_{x=1}^{N} \exp\big((ax+b_1g_1^x+\cdots +b_rg_r^x)/p \big)\,,
\end{equation*}
where $p$ is prime and $g_1,\dots, g_r$ are primitive roots $\pmod p$. An
almost forty years old problem of L. J. Mordell asks to find a nontrivial
estimate of $S$ when at least two of the coefficients $b_1,\dots,b_r$ are not
divisible by $p$. Here we obtain a nontrivial bound of the average of these
sums when $g_1$ runs over all primitive roots $\pmod p$.
\end{abstract}

\section{Introduction}

Let $p$ be a prime number, $1\le N\le p-1$, $r$ a positive integer and consider
the exponential sum
\begin{equation}\label{eq1.1}
 S_N(a,{\bb, \gg}):=
\sum_{x=1}^{N} e_p\big(ax+b_1g_1^x+\cdots +b_rg_r^x \big)\,,
\end{equation}
where $a$, and the components of $\bb=(b_1,\dots,b_r)$ are integers,
$b_1,\dots,b_r$ are not divisible by $p$ and
$\gg=(g_1, \dots,g_r)$ has components primitive roots modulo $p$. 
(We use I. M. Vinogradovs's notation $e_p(\alpha):=\exp(2\pi i \alpha/p)$.)
In the case $r=1$, when $p\mid a$ and $p\nmid b$, R. G. Stoneham\cite{Stoneham}
proved that
\begin{equation}\label{eq1.2}
 S_N(b, g):=
\sum_{x=1}^{N} e_p\big(bg^x\big)=O(p^{1/2}\log p)\,.
\end{equation}
In a correspondence with D. A. Burgess, L. J. Mordell was informed that both
Stoneham and Burgess have found 
independently several proofs of \eqref{eq1.2}. Mordell \cite{Mordell}
rediscovered one of the proofs of Burgess and observed that this leads to the
following generalization:
\begin{equation}\label{eq1.3}
 S_N(a,b, g):=
\sum_{x=1}^{N} e_p\big(ax+bg^x\big)<2p^{1/2}\log p+2p^{1/2}+1\,,
\end{equation}
where $p\nmid ab$. 
He remarks that his method doesn't seem to apply for the estimate of
\eqref{eq1.1}  when $r\ge 2$, and
the problem remained unsolved till this day. In this paper, fixing all but one
of the primitive roots, say $g\in\{g_1,\dots,g_r\}$, we derive a nontrivial
bound of $S= S_N(a,{\bb, \gg})$
on average over all $g$ primitive roots $\pmod p$.

In the following we write shortly $\gg^x=(g_1^x,\dots,g_r^x)$, for any integer
$x$ and $\gg=(g_1,\dots,g_r)$. Also,  we use the dot product notation:
$\bb\gg^x=b_1g_1^x+\cdots +b_rg_r^x$,  where $\bb=(b_1,\dots,b_r)$.
Let
\begin{equation}\label{defS}
S_N(a,b,\bb,g,\gg)=\frac{1}{\varphi(p-1)}
\,\,\sideset{}{'}\sum_{g\!\!\!\!\!\pmod p}\
\sum_{x=1}^{N} e_p\big(ax+bg^x+\bb\gg^x \big)\,,
\end{equation}
where the prime indicates that the summation is over all $g$ primitive roots
$\pmod p$.

\begin{thm}\label{Theorem1}
Let $p$ be prime, $1\le N\le p-1$, let $a,b,b_1,\dots,b_r $ be integers not all
divisible by $p$, $\gcd(b,p)=1$, and let $g,g_1,\dots,g_r$
be primitive roots $\pmod p$. 
Then:
\begin{equation}\label{eqmain}
   \big|S_N(a,b,\bb,g,\gg)\big|
\ll p^{\frac {23}{24}+\epsilon}\,.
\end{equation}
\end{thm}

The idea of proof is inspired from the Vinogradov's method and it proved
successfully in the estimation of some exponential function analogue
of Kloosterman sum, Shparlinski~\cite{Shparlinski}.
 
\section{The Complete Interval Case}
We may assume that  $r\ge 1$, since otherwise \eqref{eq1.3} gives a better
estimate than \eqref{eqmain}.
Taking some fixed primitive root $g_0 \mod p$, then any primitive root
$g$ $\pmod p$ can be written as $g=g_0^u\pmod p$, for some $1\le u\le p-1$ with
$\gcd(u,p-1)=1$. This allows us to replace the sum over $g$ in \eqref{defS} by a
sum over $1\le u\le p-1$ with $\gcd(u,p-1)=1$.
Then
\begin{equation}\label{defSS}
\begin{split}
S_N(a,b,\bb,g,\gg)=& \frac{1}{\varphi(p-1)}
\sum_{\substack{u=1\\ \gcd(u,p-1)=1}}^{p-1}
\sum_{x=1}^{N} e_p\big(ax+bg^{ux}+\bb\gg^x \big)\,,\\
\ll&\frac{\Sigma_N}{\varphi(p-1)}\,,
\end{split}
\end{equation}
where
\begin{equation*}
\Sigma_N=\sum_{\substack{u=1\\ \gcd(u,p-1)=1}}^{p-1}
\Bigg|\sum_{x=1}^N e_p\big(ax+bg^{ux}+\bb\gg^x\big)\Bigg|\,.
\end{equation*}
From now on, in this section we assume that $N=p-1$ and write shortly
$\Sigma=\Sigma_{p-1}$.
Applying the Cauchy-Schwarz inequality, we have:
\begin{equation*}
\begin{split}
\Sigma^2&\le \varphi(p-1)\sum_{\substack{u=1\\ \gcd(u,p-1)=1}}^{p-1}
\Bigg|\sum_{x=1}^{p-1} e_p\big(ax+bg^{ux}+\bb\gg^x\big)\Bigg|^2\\
&=\varphi(p-1)\sum_{\substack{u=1\\ \gcd(u,p-1)=1}}^{p-1}
\sum_{x=1}^{p-1}\sum_{y=1}^{p-1}
e_p\big(ax+\bb\gg^x-ay-\bb\gg^y\big)
e_p\big(bg^{ux}-bg^{uy}\big)\\
&\le\varphi(p-1)\sum_{x=1}^{p-1}\sum_{y=1}^{p-1}
\bigg|e_p\big(a(x-y)+\bb\gg^x-\bb\gg^y\big)\bigg|\cdot
\Bigg|\sum_{\substack{u=1\\ \gcd(u,p-1)=1}}^{p-1}e_p\big(b(g^{ux}-g^{uy})\big)
\Bigg|\,.
\end{split}
\end{equation*}
Then, by the H\"older Inequality, we get
\begin{equation*}
\begin{split}
\Sigma^8&\le \varphi(p-1)^4\Bigg(
\sum_{x=1}^{p-1}\sum_{y=1}^{p-1}
\bigg|\sum_{\substack{u=1\\ \gcd(u,p-1)=1}}^{p-1}e_p\big(b(g^{ux}-g^{uy})\big)
\bigg|\Bigg)^4\\
&\le \varphi(p-1)^4\Bigg(\sum_{x=1}^{p-1}\sum_{y=1}^{p-1}1\Bigg)^3
\sum_{x=1}^{p-1}\sum_{y=1}^{p-1}
\bigg|\sum_{\substack{u=1\\ \gcd(u,p-1)=1}}^{p-1}e_p\big(b(g^{ux}-g^{uy})\big)
\bigg|^4\,.
\end{split}
\end{equation*}
Replacing $y$ by $xy$ and then $g^x$ by $\lambda$, we have:
\begin{equation}\label{eqS8}
\begin{split}
\Sigma^8&\le p^{10}
\sum_{x=1}^{p-1}\sum_{y=1}^{p-1}
\bigg|\sum_{\substack{u=1\\ \gcd(u,p-1)=1}}^{p-1}e_p\big(b(g^{ux}-g^{uxy})\big)
\bigg|^4\\
&\le p^{10}
\sum_{\lambda=1}^{p-1}\sum_{y=1}^{p-1}
\bigg|\sum_{\substack{u=1\\
\gcd(u,p-1)=1}}^{p-1}e_p\big(b(\lambda^u-\lambda^{uy})\big)
\bigg|^4\,.
\end{split}
\end{equation}

The double sum on $y$ and $u$ can be estimated following the proof of Theorem 8
from Canetti et all~\cite{CFKLLS}. The result is
stated in the following lemma:

\begin{lem}\label{Lemma1}
For any integers $b$, $\gcd(a,b,p)=1$ and $\lambda$ primitive root $\mod p$, we
have
\begin{equation}\label{eq2.1}
H_{a,b}=\sum_{y=1}^{p-1}\Bigg|\sum_{\substack{x=1\\
\gcd(x,p-1)=1}}^{p-1}e_p\big(a\lambda^x+b\lambda^{xy}\big)\Bigg|^4
=
O\big(p^{14/3+\epsilon}\big)\,.
\end{equation}

\end{lem}
The estimate \eqref{eq2.1} is a generalization and improvement of Theorem 10
from Canetti, Friedlander, Shparlinski~\cite{CFS}. 
\begin{proof}
Using the properties of the M\"obius function and then the H\"older inequality,
we have:
\begin{equation}\label{eqH}
   \begin{split}
  H_{a,b}
=&\sum_{y=1}^{p-1}\Bigg|\sum_{d\mid p-1}\mu(d)\sum_{\substack{x=1\\
d\mid x}}^{p-1}e_p\big(a\lambda^x+b\lambda^{xy}\big)\Bigg|^4\\
\le&
\sum_{y=1}^{p-1}\Bigg(\sum_{d\mid p-1}\Bigg|
\sum_{x=1}^{\frac{p-1}{d}}e_p\big(a\lambda^{dx}+b\lambda^{dxy}\big)\
\Bigg|\Bigg)^4\\
\le&\sum_{y=1}^{p-1}\Bigg(\sum_{d\mid p-1}1\Bigg)^3
\sum_{d\mid p-1}\Bigg|
\sum_{x=1}^{\frac{p-1}{d}}e_p\big(a\lambda^{dx}+b\lambda^{dxy}\big)
\Bigg|^4\\
\le&\sigma^3_0(p-1)\sum_{d\mid p-1}\sum_{y=1}^{p-1}
\Bigg|
\sum_{x=1}^{\frac{p-1}{d}}e_p\big(a\lambda^{dx}+b\lambda^{dxy}\big)
\Bigg|^4,
   \end{split}
\end{equation}
where  $\sigma_{r}(n)=\sum_{d\mid n}d^r$  is the sum
of powers of the divisors of $n$.
For any $d\mid p-1$ we denote $t_d:=(p-1)/d$ and $\lambda_d:=\lambda^d$.
Notice that the multiplicative order of $\lambda^d$ is $t_d$. Then  the sum
over $y$ from the last line above becomes:
\begin{equation}\label{eqH2}
\begin{split}
\sum_{y=1}^{p-1}
\Bigg|
\sum_{x=1}^{\frac{p-1}{d}}e_p\big(a\lambda^{dx}+b\lambda^{dxy}\big)
\Bigg|^4
=&d\sum_{y=1}^{t_d}\Bigg|\sum_{x=1}^{t_d}e_p\big(a\lambda_d^{x}+b\lambda_d^{xy}
\big)\Bigg|^4\\
=&d\sum_{y=1}^{t_d}\frac{1}{t_d}\sum_{z=1}^{t_d}
\Bigg|\sum_{x=1}^{t_d}
e_p\big(a\lambda_d^{x+z}+b\lambda_d^{(x+z)y} \big)
\Bigg|^4\\
=&\frac{d}{t_d}\sum_{y=1}^{t_d}\sum_{z=1}^{t_d}
\Bigg|\sum_{x=1}^{t_d}
e_p\big(a\lambda_d^{z}\lambda_d^{x}+b\lambda_d^{zy}\lambda_d^{xy} \big)
\Bigg|^4\\
\le&\frac{d}{t_d}\sum_{y=1}^{t_d}\sum_{\alpha,\beta=0}^{p-1}
\Bigg|\sum_{x=1}^{t_d}
e_p\big(\alpha\lambda_d^{x}+\beta\lambda_d^{xy} \big)
\Bigg|^4,
\end{split}
\end{equation}
since for each fixed $y\in\{1,\dots,t_d\}$ the pairs
$(at_d^z,b\lambda_d^{zy})$ with $z\in\{1,\dots,t_d\}$ are distinct modulo $p$.
Next we write explicitly the absolute value in the last term and see that
\eqref{eqH2} gives
\begin{equation}\label{eqH3}
\begin{split}
&\sum_{y=1}^{p-1}
\Bigg|
\sum_{x=1}^{\frac{p-1}{d}}e_p\big(a\lambda^{dx}+b\lambda^{dxy}\big)
\Bigg|^4\\
\le&\frac{d}{t_d}
\sum_{\alpha,\beta=0}^{p-1}\
\sum_{y,x_1,x_2,x_3,x_4=1}^{t_d}
e_p\big(\alpha(\lambda_d^{x_1}+\lambda_d^{x_2}-\lambda_d^{x_3}-\lambda_d^{x_4})
+\beta(\lambda_d^{ x_1y}+\lambda_d^{ x_2y}-\lambda_d^{ x_3y}-\lambda_d^{ x_4y})
\big)
\\
\le&\frac{d}{t_d}\ 
\sum_{\substack{y,x_1,x_2=1\\x_3,x_4=1}}^{t_d}
\sum_{\alpha=0}^{p-1}
e_p\big(\alpha(\lambda_d^{x_1}+\lambda_d^{x_2}-\lambda_d^{x_3}-\lambda_d^{x_4}
)\big)
\sum_{\beta=0}^{p-1}
e_p\big(\beta(\lambda_d^{ x_1y}+\lambda_d^{ x_2y}-\lambda_d^{ x_3y}-\lambda_d^{
x_4y})
\big)
\\
=&\frac{d}{t_d}\cdot p^2\cdot T_d\,,
\end{split}
\end{equation}
where $T_d$ is the number of solutions of the system of congruences:
\begin{equation*}
   \begin{cases}
      \lambda_d^{x_1}+\lambda_d^{x_2}\equiv \lambda_d^{x_3}+\lambda_d^{x_4},\\
       \lambda_d^{x_1y}+\lambda_d^{x_2y}\equiv
\lambda_d^{x_3y}+\lambda_d^{x_4y}\\	
   \end{cases}
\end{equation*}
with $1\le x_1,x_2,x_3,x_4,y\le t_d$.
In the proof of Theorem 8 from Canetti et all~\cite{CFKLLS}, 
the last inequality bounds $T_d$ by
\begin{equation}\label{eqTd}
   T_d\ll t_d^{14/3}p^{-1}\,.
\end{equation}
Then, by \eqref{eqH},  \eqref{eqH3} and \eqref{eqTd}, we obtain
\begin{equation*}
\begin{split}
   H_{a,b}\ll &\ \sigma^3_0(p-1)\sum_{d\mid p-1}
\frac {d}{t_d}\cdot p^2\cdot t_d^{14/3}p^{-1}\\
\ll& \ \sigma^3_0(p-1)\sigma_{-\frac{14}{3}}(p-1)\cdot p^{14/3}\,
\end{split}
\end{equation*}
and the lemma follows, since $\sigma_r(n)\ll n^ \epsilon$ for any $r$.

\end{proof}

By \eqref{eqS8} and \eqref{eq2.1} we deduce that:
\begin{equation*}
\Sigma^8\ll p^{10}
\sum_{\lambda=1}^{p-1}p^{14/3+\epsilon}\ll p^{47/3+\epsilon}\,.
\end{equation*}
Then making use of the estimate
$p/\log\log p \ll \varphi(p-1)$, we obtain
\begin{equation}\label{eqSgata}
\frac{\Sigma}{\varphi(p-1)}\ll p^{23/24+\epsilon}\,.
\end{equation}
From this estimate together with \eqref{defSS}, it follows \eqref{eqmain}, so
Theorem~\ref{Theorem1} is proved in the case $N=p-1$.

\section{Completion of the Proof}

It remains to show that the size of the incomplete sums is not far from that of
the complete ones.
Let $I$ be an interval of integers $\subseteq [1,p-1]$ and denote 
\begin{equation}\label{SI}
 S(I)=\sum_{\substack{u=1\\ \gcd(u,p-1)=1}}^{p-1}
\sum_{x\in I}e_p(ax+bg^{ux}+\bb\gg^x)\,.
\end{equation}
In order to estimate the departure of $S(I)$ from $S([1,p-1])$, 
the following characteristic function of the interval $I$ is suitable:
\begin{equation*}
 \frac 1p \sum_{y\in I} \sum_{k=1}^{p}e_p\big(k(y-x)\big)
=\begin{cases}
1, \quad \text{if } x\in I;\\
0, \quad \text{else}.
 \end{cases}
\end{equation*}
Then
\begin{equation*}
\begin{split}
 S(I)=&\sum_{\substack{u=1\\ \gcd(u,p-1)=1}}^{p-1}
\sum_{x\in I}e_p(ax+bg^{ux}+\bb\gg^x)\\
=&\sum_{\substack{u=1\\ \gcd(u,p-1)=1}}^{p-1}
\sum_{x=1}^{p-1}e_p(ax+bg^{ux}+\bb\gg^x)
 \frac 1p \sum_{y\in I} \sum_{k=1}^{p}e_p\big(k(y-x)\big)\\
=& \frac 1p \sum_{k=1}^{p} \sum_{y\in I} e_p(ky)
\sum_{\substack{u=1\\ \gcd(u,p-1)=1}}^{p-1}
\sum_{x=1}^{p-1}e_p\big((a-k)x+bg^{ux}+\bb\gg^x\big)\,.
\end{split}
\end{equation*}
In this last form of $S(I)$ we separate the terms with $k=p$ and bound its
absolute value to get:
\begin{equation}\label{eqSIag}
\begin{split}
 |S(I)|\le&\ \frac 1p \sum_{k=1}^{p-1}\Bigg| \sum_{y\in I} e_p(ky) \Bigg|
\sum_{\substack{u=1\\ \gcd(u,p-1)=1}}^{p-1}
\Bigg|
\sum_{x=1}^{p-1}e_p\big((a-k)x+bg^{ux}+\bb\gg^x\big)
\Bigg|\\
&+\frac 1p |I| \sum_{\substack{u=1\\ \gcd(u,p-1)=1}}^{p-1}
\Bigg|
\sum_{x=1}^{p-1}e_p(ax+bg^{ux}+\bb\gg^x)
\Bigg|\,.
\end{split}
\end{equation}
Here the sum over $y$ is a geometric progression, that can be evaluated
accurately using
\begin{equation}\label{eqsinus}
 \big| e_p(k)-1 \big| =  2\Big| \sin \Big( \frac{k\pi}{p} \Big) \Big|
\ge 4\,\norm{\frac{k}{p}}\,,
\end{equation}
where $\norm{\cdot}$ is the distance to the nearest integer,
while the sums over $u$ and $x$ are the complete sums bounded by
\eqref{eqSgata}.
Thus, by \eqref{eqSIag}, \eqref{eqsinus} and \eqref{eqSgata}, we get
\begin{equation*}
\begin{split}
|S(I)|\le &\frac 1p \sum_{k=1}^{p-1}
\frac{2}{\big|e_p(k)-1\big|} p^{47/24+\epsilon}
+\frac 1p|I|p^{47/24+\epsilon}\\
\le& p^{47/24+\epsilon}
\Bigg(\frac 1p \sum_{k=1}^{\frac{p-1}{2}}\frac{1}{2k/p}+1\Bigg)\\
\le& p^{47/24+\epsilon}(3+\log p)\\
\le& p^{47/24+\epsilon}\,,
\end{split}
\end{equation*}
which concludes the proof of Theorem~\ref{Theorem1}.

\bigskip
{\bf Acknowledgement: }
We thank Ke Gong for his careful reading of an earlier version of this paper.

\end{document}